\documentclass[12pt]{amsart}
\usepackage{amssymb}
\usepackage{verbatim}
\usepackage{mathrsfs}
\usepackage{mathtools}
\usepackage{caption}

\usepackage{xcolor}

\newtheorem{thm}{Theorem}[section]
\newtheorem{cor}[thm]{Corollary}

\numberwithin{equation}{section}

\newcommand{\alt}{\alpha}
\newcommand{\Fs}{ F_2 }
\newcommand{\Lt}{ l_2 }
\newcommand{\Aa}{\mathbf{a}}
\newcommand{\HyG}{ {}_2F_1 }

\providecommand{\sym}{\operatorname{sym}}

\mathtoolsset{showonlyrefs}

\begin{document}

\title[Asymptotics of a Gauss hypergeometric function]{Asymptotics of a Gauss hypergeometric function related to moments of symmetric-square $L$-functions II}

\author[D. Frolenkov]{Dmitry  Frolenkov}
\address{
Steklov Mathematical Institute of Russian Academy of Sciences, 8 Gubkina st., Moscow, 119991, Russia}
\email{frolenkov@mi-ras.ru}

\begin{abstract}
We prove an asymptotic formula for $\HyG\left(1/4-it+ir,1/4-it-ir,1/2;x \right)$ as $r,  t\to\infty$ and $\alt=r/t\to0.$
This special case of the Gauss hypergeometric function appears in the explicit formula for the first moment of Maass form symmetric-square $L$-functions $L(\sym^2 u_{j},1/2+2it)$.
\end{abstract}

\keywords{hypergeometric function, Bessel function, Airy function, Liouville-Green method}
\subjclass[2010]{Primary: 33C05, 41A60}

\maketitle


\section{Introduction}
Consider the function
\begin{multline}\label{2f1 1/4 1/4 1/2 def}
\Fs(r,\alt,x):=
\frac{\Gamma(1/4-it+ir)\Gamma(1/4-it-ir)}{\Gamma(1/2)}\HyG\left(1/4-it+ir,1/4-it-ir,1/2;x \right)=\\=
\frac{\Gamma(1/4+ir(1-\alt))\Gamma(1/4-ir(1+\alt))}{\Gamma(1/2)}\HyG\left(1/4+ir(1-\alt),1/4-ir(1+\alt),1/2;x \right).
\end{multline}
This hypergeometric function appears in the explicit formula for moments of Maass form symmetric square $L$-functions on the critical line $\rho=1/2+2it$ ( see \cite[Theorem 5, Lemma 7]{Bal} for details) and plays a crucial role in the asymptotic analysis of the corresponding moments. In view of this application, our main result is the following.
\begin{thm}\label{thm: 2F1 asympt1}
For $0<y<1,$ $r\to\infty$ and  $r^{-1+\delta}<\alt< r^{-\delta}$ one has
\begin{multline}\label{2F1 main asymptI}
\HyG\left(1/4+ir(1-\alt),1/4-ir(1+\alt),1/2;y \right)=\frac{(1-y)^{ir\alt}}{2\sqrt{\pi}}\left(\frac{(1-\alt^2)(\zeta-\alt^2)r^2}{\zeta^2(1-\alt^2-y)}\right)^{1/4}\times\\\times
\Biggl(
2e^{\pi r}W_{2n+1,2}(2r,\zeta)\left(1+\sum_{j=1}^{n-1}\frac{e_j(\alt)}{(r\alt)^j}+O((r\alt)^{-n})\right)+\\
+e^{-\pi r+2\pi r\alt}W_{2n+1,1}(2r,\zeta)\left(1+\sum_{j=1}^{n-1}\frac{f_j(\alt)}{(r\alt)^j}+O((r\alt)^{-n})\right)
\Biggr),
\end{multline}
where
\begin{multline}\label{W1 def}
W_{2n+1,1}(\mu,\zeta)=\sqrt{\zeta}\tilde{I}_{i\mu\alt}(\mu\sqrt{\zeta})\sum_{s=0}^{n}\frac{A_s(\alt,\zeta)}{\mu^{2s}}+
\frac{\zeta}{\mu}\tilde{I}'_{i\mu\alt}(\mu\sqrt{\zeta})\sum_{s=0}^{n-1}\frac{B_s(\alt,\zeta)}{\mu^{2s}}+\epsilon_{2n+1,1}(\mu,\alt,\zeta),
\end{multline}
\begin{multline}\label{W2 def}
W_{2n+1,2}(\mu,\zeta)=\sqrt{\zeta}K_{i\mu\alt}(\mu\sqrt{\zeta})\sum_{s-0}^{n}\frac{A_s(\alt,\zeta)}{\mu^{2s}}+
\frac{\zeta}{\mu}K'_{i\mu\alt}(\mu\sqrt{\zeta})\sum_{s=0}^{n-1}\frac{B_s(\alt,\zeta)}{\mu^{2s}}+\epsilon_{2n+1,2}(\mu,\alt,\zeta),
\end{multline}
with $A_s(\alt,\zeta)$, $B_s(\alt,\zeta)$ defined by \cite[(7.19)-(7.21)]{Dun13} and
\begin{equation}
\epsilon_{2n+1,1}(\mu,\alt,\zeta)=\sqrt{\zeta}\tilde{I}_{i\mu\alt}(\mu\sqrt{\zeta})O\left(\mu^{-2n-1}\right),\quad
\epsilon_{2n+1,2}(\mu,\alt,\zeta)=\sqrt{\zeta}K_{i\mu\alt}(\mu\sqrt{\zeta})O\left(\mu^{-2n-1}\right).
\end{equation}
The coefficients $e_j(\alt),f_j(\alt)\ll1$ can be expressed in terms of $A_s, B_s$. For $0<y<1-\alt^2$
\begin{equation}\label{zeta>alt^2}
\arccos\frac{\sqrt{y}}{\sqrt{1-\alt^2}}-\frac{\alt}{2}\arccos\left(\frac{y(1+\alt^2)-(1-\alt^2)}{(1-\alt^2)(1-y)}\right)=
\sqrt{\zeta-\alt^2}-\alt\arccos\frac{\alt}{\sqrt{\zeta}},
\end{equation}
for  $1-\alt^2\le y<1$
\begin{multline}\label{zeta<alt^2}
\log\frac{\sqrt{y}-\sqrt{y-1+\alt^2}}{\sqrt{1-\alt^2}}-
\frac{\alt}{2}\log\left(\frac{2\alt\left(\alt-\sqrt{y(y-1+\alt^2)}\right)}{(1-\alt^2)(1-y)}-\frac{1+\alt^2}{1-\alt^2}\right)=\\=
\alt\log\left(\frac{\alt+\sqrt{\alt^2-\zeta}}{\sqrt{\zeta}}\right)-\sqrt{\alt^2-\zeta},
\end{multline}
and $\tilde{I}_{i\nu}(x)=\pi e^{-\pi\nu}\left(I_{i\nu}(x)+I_{-i\nu}(x)\right)$.
\end{thm}
To simplify  further application of this result, we formulate several corollaries.

\begin{cor}\label{cor: y<1-alt^2}
For $0<y<1-\alt^2-\frac{\alt^{4/3}}{r^{2/3-\delta}}$ the estimate takes place:
\begin{equation}\label{2f1 1/4 1/4 1/2 asympt y<1-alt^2}
\Fs(r,\alt,y)\ll
\frac{e^{-\pi r\alt}e^{-2r\Aa_0(\alt,y)}}{\sqrt{r}(1-\alt^2-y)^{1/4}}\ll\frac{e^{-\pi r\alt}}{r^{A}},
\end{equation}
where
\begin{equation}\label{Aa0 def}
\Aa_0(\alt,y)=
\frac{\pi(1-\alt)}{2}-\arctan\frac{\sqrt{y}}{\sqrt{1-\alt^2-y}}+\alt\arctan\frac{\alt\sqrt{y}}{\sqrt{1-\alt^2-y}}.
\end{equation}
\end{cor}
\begin{cor}\label{cor: y>1-alt^2}
For $1-\alt^2+\frac{\alt^{4/3}}{r^{2/3-\delta}}<y<1$ the leading coefficient in the asymptotic expansion for $\Fs(r,\alt,y)$ is given by
\begin{equation}\label{2f1 1/4 1/4 1/2 asympt y>1-alt^2}
\Fs(r,\alt,y)\sim
2\sqrt{\pi}e^{ir\Lt(\alt,y)}e^{-\pi r\alt}
\frac{\cos\left(2r\Aa_1(\alt,y)-\pi/4\right)}{\sqrt{r}(y-1+\alt^2)^{1/4}},
\end{equation}
where
\begin{equation}\label{Aa1 def}
\Aa_1(\alt,y)=
\alt\log\left(\alt\sqrt{y}+\sqrt{y-1+\alt^2}\right)-\log\left(\sqrt{y}+\sqrt{y-1+\alt^2}\right)-\frac{\alt}{2}\log(1-y)+
\frac{1-\alt}{2}\log(1-\alt^2),
\end{equation}
\begin{equation}\label{2F1 main asymptII Lt func}
\Lt(\alt,y)=\alt\log(1-y)-2\alt\log r+(1-\alt)\log(1-\alt)-(1+\alt)\log(1+\alt)+2\alt.
\end{equation}
\end{cor}

\begin{cor}\label{cor: y=1-alt^2}
For $0<y<1$
one has
\begin{equation}\label{2f1 1/4 1/4 1/2 asympt y=1-alt^2}
\Fs(r,\alt,y)\sim
e^{ir\Lt(\alt,y)}
\frac{ 2^{3/2}\pi e^{-\pi r\alt}}{(2r\alt)^{1/3}}\left(\frac{\alt^2\hat{\zeta}(y)}{y-1+\alt^2}\right)^{1/4}Ai(-(2r\alt)^{2/3}\hat{\zeta}(y)),
\end{equation}
where
\begin{equation}\label{zeta(y) def}
\hat{\zeta}(y)=\left\{
            \begin{array}{ll}
              -\left(3\Aa_0(\alt,y)/(2\alt)\right)^{2/3}, & \hbox{if} \quad y<1-\alt^2 \\
              \left(3\Aa_1(\alt,y)/(2\alt)\right)^{2/3}, & \hbox{if} \quad y>1-\alt^2.
            \end{array}
          \right.
\end{equation}
\end{cor}


In the case when $\alt=0$ one can apply \cite[Theorem 3.2]{KD2} in order to obtain an asymptotic expansion for $\Fs(r,\alt,x)$ in terms of $K$-Bessel functions. This shows that the function \eqref{2f1 1/4 1/4 1/2 def} is exponentially small for $1-x\gg r^{-2+\epsilon}$ (see \cite[Corollary 4.2, Lemma 4.3]{BF2mom}). It is also possible to generalize \cite[Theorem 3.2]{KD2} to the case $\alt\neq0.$  Doing so, we obtain a version of \cite[(3.40)]{KD2} for \eqref{2f1 1/4 1/4 1/2 def} (thus we set $z=1-2x$ in \cite[(3.40)]{KD2}) with
\begin{equation}
f(t)=\log\left(1+\frac{1-x}{t}\right)-(1+\alt)\log(1+t)+\alt\log\left(t(t+1-x)\right).
\end{equation}
Saddle points of this function are
\begin{equation}
t_{\pm}=\frac{x-(1-\alt)\pm\sqrt{x(x-1+\alt^2)}}{1-\alt}.
\end{equation}
For $x>1-\alt^2$ these points are real, while for $0<x<1-\alt^2$ they are not. Therefore, the function $\Fs(r,\alt,x)$ has different asymptotic behaviour in the regions $x<1-\alt^2$ and $x>1-\alt^2$. In fact, the shape of $f(t)$ is quite different for $x\lessgtr1-\alt^2.$ If $x>1-\alt^2$ it resembles a plot of a cubic  polynomial, which suggests that $\Fs(r,\alt,x)$ may be approximated in terms of the Airy function. While for  $x<1-\alt^2$ the shape of $f(t)$  suggests approximation in terms of Bessel functions. Usually, when shapes of two asymptotic expansions to the left and to the right of some point are different, both of them are unlikely to be uniformly valid in a neighbourhood of this point. Therefore, it is required to construct a somewhat more complex asymptotic expansion that is valid for all $x$ near  $1-\alt^2.$

In order to obtain such a uniform asymptotic expansion, we decided to apply the Liouville-Green method.  It follows from \cite[p.96, (7)-(9)]{BE} that
the function
\begin{equation}\label{Y def}
Y_r(x)=x^{1/4}(1-x)^{1/2-ir\alt}\HyG\left(1/4+ir(1-\alt),1/4-ir(1+\alt),1/2;x \right)
\end{equation}
satisfies the differential equation
\begin{equation}\label{Y difeq}
Y_r''(x)=\left((2r)^2f(\alt,x)+g(\alt,x)\right)Y_r(x),
\end{equation}
\begin{equation}\label{Y difeq fg def}
f(\alt,x)=\frac{1-\alt^2-x}{4x(1-x)^2},\quad
g(\alt,x)=-\frac{1}{4(1-x)^2}-\frac{3}{16x^2(1-x)}.
\end{equation}
This differential equation has a double pole at $x=1$ and a turning point at $x=1-\alt^2.$ These points coalesce  as $\alt\to0$. Asymptotic properties of real solutions of such differential equations have been studied by Boyd-Dunster \cite{BD} and  Dunster \cite{Dun90}.  Note that as $x\to1$ one has
\begin{equation}
(1-x)^2f(\alt,x)\to-\alt^2/4,\quad (1-x)^2g(\alt,x)\to-1/4,
\end{equation}
and therefore, the differential equation \eqref{Y difeq} is exactly as the one studied in \cite{Dun90} (see conditions on \cite[p.1595 lines 1-5]{Dun90}).
We remark that the function $Y_r(x)$ is real-valued since
\begin{equation}
Y_r(x)=x^{1/4}(1-x)^{1/2+ir\alt}\HyG\left(1/4-ir(1-\alt),1/4+ir(1+\alt),1/2;x \right)=Y_{-r}(x),
\end{equation}
which follows from \cite[15.8.1]{HMF}.

In \cite[Theorem 1]{Dun90}  asymptotic expansions for solutions of \eqref{Y difeq}  are given in terms of $K_{2ir\alt}(\cdot)$ and $\tilde{I}_{2ir\alt}(\cdot)$ Bessel functions. Thus it is left to identify  $Y_r(x)$ with one of these functions. The problem is that there does not exist a point $x_0$ where both $Y_r(x)$ and one of the functions $K_{2ir\alt}(\cdot)$ or $\tilde{I}_{2ir\alt}(\cdot)$ are recessive. This means that $Y_r(x)$ is not a multiple of one of the solutions $W_{2n+1,j}$ from \cite[Theorem 1]{Dun90}, but is their linear combination. This is not surprising since in the case when $\alt=0$ (see \cite[ Lemma 4.1, Corollary 4.2]{BF2mom}) the asymptotic expansion for $Y_r(x)$ is given in terms of $K_0(\cdot)$ and $I_0(\cdot)$ Bessel functions. There are several methods to find coefficients in this linear combination expansion, but all of them are quite tricky. For example, one can argue in the same way as in  \cite[sec. 5.1]{BFJEMS}. Luckily, this issue was addressed in \cite[sec.7]{Dun13}, where a function similar to $Y_r(x)$ has been considered, namely \cite[(1.9)]{Dun13}
\begin{equation}\label{e def}
e^{i\mu}_{-1/2+i\tau}(x)=(1+x^2)^{-i\mu/2}\HyG\left(\frac{1}{4}-\frac{i(\mu+\tau)}{2},\frac{1}{4}-\frac{i(\mu-\tau)}{2},1/2;-x^2 \right).
\end{equation}
The function above is closely related to \eqref{2f1 1/4 1/4 1/2 def}, since
applying \cite[15.8.1]{HMF} one has
\begin{equation}\label{2F1 to e}
\HyG\left(1/4+ir(1-\alt),1/4-ir(1+\alt),1/2;y \right)=
(1-y)^{-1/4+ir\alt}e^{2ir}_{-1/2+2ir\alt}\left(\sqrt{\frac{y}{1-y}}\right).
\end{equation}
Though the function $e^{i\mu}_{-1/2+i\tau}(x)$ with $\mu\to\infty$ and $\tau=\mu\alt<(1-\delta)\mu$ has been studied in \cite[sec.7]{Dun13} and almost all required computations were performed,  the final asymptotic expansion  is not written down in \cite[sec.7]{Dun13}.  We are going to finalize computations in  \cite[sec.7]{Dun13} keeping in mind that we are looking for an asymptotic expansion of the left-hand side of \eqref{2F1 to e} (we will formulate some results both in $x$ and $y$ notations eventually substituting $x=y^{1/2}(1-y)^{-1/2}$).

\section{Some results of Dunster}
In order to express  $e^{i\mu}_{-1/2+i\tau}(x)$ as a linear combination of solutions constructed in \cite[Theorem 1]{Dun90} it is first required to introduce a second function satisfying the same differential equation and completing $e^{i\mu}_{-1/2+i\tau}(x)$ to a numerically satisfactory pair of solutions. To this end, Dunster introduced
\cite[(1.10)]{Dun13}
\begin{equation}\label{0 def}
o^{i\mu}_{-1/2+i\tau}(x):=x(1+x^2)^{-i\mu/2}\HyG\left(\frac{3}{4}-\frac{i(\mu+\tau)}{2},\frac{3}{4}-\frac{i(\mu-\tau)}{2},3/2;-x^2 \right).
\end{equation}
Performing the change of variables:
\begin{equation}\label{sxy def}
s:=1-\frac{x}{\sqrt{x^2+1}}=1-\sqrt{y},
\end{equation}
and defining $V(s):=(1+x^2)^{-3/4}\omega(x)$, where $\omega(x)$ is either $e^{i\mu}_{-1/2+i\tau}(x)$ or $o^{i\mu}_{-1/2+i\tau}(x)$, one has (see \cite[(7.2)]{Dun13})
\begin{equation}\label{V difeq}
\frac{d^2V}{ds^2}=\left(\frac{\mu^2(2s-s^2-\alt^2)}{s^2(2-s)^2}+\frac{2s-s^2-4}{4s^2(2-s)^2}\right)V.
\end{equation}
Applying the theory developed in \cite{Dun90}, Dunster obtained the system \cite[(7.30)]{Dun13} of two linear equations for $e^{i\mu}_{-1/2+i\tau}(x)$ and $o^{i\mu}_{-1/2+i\tau}(x)$. Solving it one has
\begin{equation}\label{e to W1W2}
e^{i\mu}_{-1/2+i\tau}(x)=\left(\frac{\zeta-\alt^2}{\zeta^2(1-\alt^2-\alt^2x^2)}\right)^{1/4}
\frac{W_{2n+1,1}(\mu,\zeta)C^o_{2n+1,2}(\mu)-W_{2n+1,2}(\mu,\zeta)C^o_{2n+1,1}(\mu)}{C^e_{2n+1,1}(\mu)C^o_{2n+1,2}(\mu)-C^e_{2n+1,2}(\mu)C^o_{2n+1,1}(\mu)},
\end{equation}
where $\zeta$ is a new variable defined for $0<x\le(1-\alt^2)^{1/2}\alt^{-1}$ (or equivalently for $0<y<1-\alt^2$) by \cite[(7.6)]{Dun13}
\begin{multline}\label{zeta>alt^2}
\arccos\frac{x}{\sqrt{(1-\alt^2)(1+x^2)}}-\frac{\alt}{2}\arccos\left(\frac{2\alt^2(1+x^2)-(1+\alt^2)}{1-\alt^2}\right)=\\=
\arccos\frac{\sqrt{y}}{\sqrt{1-\alt^2}}-\frac{\alt}{2}\arccos\left(\frac{y(1+\alt^2)-(1-\alt^2)}{(1-\alt^2)(1-y)}\right)=
\sqrt{\zeta-\alt^2}-\alt\arctan\frac{\sqrt{\zeta-\alt^2}}{\alt}=\\=
\sqrt{\zeta-\alt^2}-\alt\arccos\frac{\alt}{\sqrt{\zeta}}
\end{multline}
and for $x\ge(1-\alt^2)^{1/2}\alt^{-1}$ (or equivalently for $1-\alt^2\le y<1$) by \cite[(7.7)]{Dun13}
\begin{multline}\label{zeta<alt^2}
\log\left(\frac{x-\sqrt{1+\alt^2(1+x^2)}}{\sqrt{(x^2+1)(1-\alt^2)}}\right)-
\frac{\alt}{2}\log\left(\frac{2\alt}{1-\alt^2}\left(\alt(x^2+1)-x\sqrt{1+\alt^2(1+x^2)}\right)-\frac{1+\alt^2}{1-\alt^2}\right)=\\=
\log\frac{\sqrt{y}-\sqrt{y-1+\alt^2}}{\sqrt{1-\alt^2}}-
\frac{\alt}{2}\log\left(\frac{2\alt\left(\alt-\sqrt{y(y-1+\alt^2)}\right)}{(1-\alt^2)(1-y)}-\frac{1+\alt^2}{1-\alt^2}\right)=\\
=\frac{\alt}{2}\log\left(\frac{\alt+\sqrt{\alt^2-\zeta}}{\alt-\sqrt{\alt^2-\zeta}}\right)-\sqrt{\alt^2-\zeta}=
\alt\log\left(\frac{\alt+\sqrt{\alt^2-\zeta}}{\sqrt{\zeta}}\right)-\sqrt{\alt^2-\zeta},
\end{multline}
where $W_{2n+1,j}(\mu,\zeta)$ are defined by \eqref{W1 def}, \eqref{W2 def} (see \cite[(7.16),(7.17)]{Dun13})
with coefficients $A_s(\alt,\zeta)$ and $B_s(\alt,\zeta)$ given by \cite[(7.19)-(7.21)]{Dun13}. Various properties of the error terms $\epsilon_{2n+1,j}(\mu,\alt,\zeta)$ are stated in \cite[(7.22)-(7.25)]{Dun13}. Finally, according to \cite[(7.31), (7.32)]{Dun13} one has
\begin{equation}\label{Ce def}
C^e_{2n+1,j}(\mu)=\left(\frac{\zeta_0-\alt^2}{\zeta_0^2(1-\alt^2)}\right)^{1/4}W_{2n+1,j}(\mu,\zeta_0),
\end{equation}
\begin{equation}\label{Co def}
C^o_{2n+1,j}(\mu)=
\frac{-2\zeta_0(\zeta_0-\alt^2)W_{2n+1,j}'(\mu,\zeta_0)+(\zeta_0/2-\alt^2)W_{2n+1,j}(\mu,\zeta_0)}{\zeta_0^{1/2}(\zeta_0-\alt^2)^{5/4}(1-\alt^2)^{3/4}},
\end{equation}
where $W_{2n+1,j}'(\mu,\zeta_0)=\frac{d}{d\zeta}W_{2n+1,j}(\mu,\zeta)\Bigl|_{\zeta=\zeta_0}$, and $\zeta_0$ is the point corresponding to $x=0$, namely (see \cite[(7.8)]{Dun13})
\begin{equation}\label{zeta_0}
\frac{\pi(1-\alt)}{2}=\sqrt{\zeta_0-\alt^2}-\alt\arctan\frac{\sqrt{\zeta_0-\alt^2}}{\alt}.
\end{equation}
\section{Proof of Theorem \ref{thm: 2F1 asympt1}}
Using \eqref{Ce def}, \eqref{Co def} we find that
\begin{multline}\label{CeCo-Ceco}
C^e_{2n+1,1}(\mu)C^o_{2n+1,2}(\mu)-C^e_{2n+1,2}(\mu)C^o_{2n+1,1}(\mu)=\\=
\frac{-2}{1-\alt^2}\left(W_{2n+1,1}(\mu,\zeta_0)W_{2n+1,2}'(\mu,\zeta_0)-W_{2n+1,1}'(\mu,\zeta_0)W_{2n+1,2}(\mu,\zeta_0)\right)=\\=
\frac{-2}{1-\alt^2}\mathcal{W}\left(W_{2n+1,1},W_{2n+1,2}\right)(\mu,\zeta_0),
\end{multline}
where $\mathcal{W}(f,g)=fg'-f'g$ is a Wronskian of two functions $f(x)$ and $g(x)$.
Let
\begin{equation}\label{Sigmaab def}
\Sigma_a(\zeta):=\sum_{s=0}^{n}\frac{A_s(\alt,\zeta)}{\mu^{2s}},\quad
\Sigma_b(\zeta):=\sum_{s=0}^{n-1}\frac{B_s(\alt,\zeta)}{\mu^{2s}},
\end{equation}
\begin{equation}\label{SigmaA def}
\acute{\Sigma}_A(\zeta):=\frac{d}{d\zeta}\left(\sqrt{\zeta}\Sigma_a(\zeta)\right)+
\left(1-\frac{\tau^2}{\mu^2\zeta}\right)\frac{\sqrt{\zeta}}{2}\Sigma_b(\zeta),
\end{equation}
\begin{equation}\label{SigmaB def}
\acute{\Sigma}_B(\zeta):=\frac{d}{d\zeta}\left(\frac{\zeta}{\mu}\Sigma_b(\zeta)\right)+
\frac{\mu}{2}\Sigma_a(\zeta)-\frac{1}{2\mu}\Sigma_b(\zeta).
\end{equation}
Then using \cite[10.25.1]{HMF}, \cite[(7.25)]{Dun13} we conclude that
\begin{equation}\label{dW2/dzeta}
\frac{d}{d\zeta}W_{2n+1,2}(\mu,\zeta)\Bigl|_{\zeta=\zeta_0}=
K_{i\tau}(\mu\sqrt{\zeta_0})\acute{\Sigma}_A(\zeta)+
\frac{d}{d\zeta}K_{i\tau}(\mu,\zeta)\Bigl|_{\zeta=\zeta_0}\acute{\Sigma}_B(\zeta)
\end{equation}
and
\begin{equation}\label{dW1/dzeta}
\frac{d}{d\zeta}W_{2n+1,1}(\mu,\zeta)\Bigl|_{\zeta=\zeta_0}=
\tilde{I}_{i\tau}(\mu\sqrt{\zeta_0})\acute{\Sigma}_A(\zeta)+
\frac{d}{d\zeta}\tilde{I}_{i\tau}(\mu,\zeta)\Bigl|_{\zeta=\zeta_0}\acute{\Sigma}_B(\zeta)+
\frac{d}{d\zeta}\epsilon_{2n+1,1}(\mu,\alt,\zeta)\Bigl|_{\zeta=\zeta_0}.
\end{equation}
Therefore,
\begin{multline}\label{Wronskian W1W2}
\mathcal{W}\left(W_{2n+1,1},W_{2n+1,2}\right)(\mu,\zeta_0)=
\mathcal{W}\left(\epsilon_{2n+1,1},W_{2n+1,2}\right)(\mu,\zeta_0)+\\+
\left(\frac{\zeta_0}{\mu}\Sigma_b(\zeta_0)\acute{\Sigma}_A(\zeta_0)-\sqrt{\zeta_0}\Sigma_a(\zeta_0)\acute{\Sigma}_B(\zeta_0)\right)
\mathcal{W}\left(K_{i\tau},\tilde{I}_{i\tau}\right)(\mu\sqrt{\zeta_0}).
\end{multline}
It follows from \cite[(2.10)]{Dun90Bessel}, \cite[(2.8)]{Dun90} that
\begin{equation}
\mathcal{W}\left(K_{i\tau},\tilde{I}_{i\tau}\right)(z)=\frac{2\pi e^{-\pi\tau}}{z}.
\end{equation}
Substituting  \eqref{SigmaA def}, \eqref{SigmaB def} to \eqref{Wronskian W1W2} and choosing $A_0=1$ we prove that
\begin{multline}\label{Wronskian W1W2 asympt}
\mathcal{W}\left(W_{2n+1,1},W_{2n+1,2}\right)(\mu,\zeta_0)=
\pi e^{-\pi\tau}\left(
-\Sigma_a^2(\zeta_0)-\frac{2\zeta_0}{\mu^2}\mathcal{W}\left(\Sigma_a,\Sigma_b\right)(\zeta_0)+
\left(1-\frac{\tau^2}{\mu^2\zeta_0}\right)\frac{\zeta_0}{\mu^2}\Sigma_b^2(\zeta_0)
\right)+\\+\mathcal{W}\left(\epsilon_{2n+1,1},W_{2n+1,2}\right)(\mu,\zeta_0)=
-\pi e^{-\pi\tau}\left(1+\sum_{s=1}^{2n}\frac{C_s(\alt,\zeta_0)}{\mu^{2s}}\right)+
\mathcal{W}\left(\epsilon_{2n+1,1},W_{2n+1,2}\right)(\mu,\zeta_0),
\end{multline}
where coefficients $C_s(\alt,\zeta)$ can be written in terms of $A_s(\alt,\zeta)$ and $B_s(\alt,\zeta)$. According to \cite[(7.22)]{Dun13} (see also \cite[(2.29),(2.36)]{Dun90})
\begin{equation}
\epsilon_{2n+1,1}(\mu,\alt,\zeta_0)\ll\frac{\left|\tilde{I}_{i\tau}(\mu\sqrt{\zeta_0})\right|}{\mu^{2n+1}},\quad
\frac{d}{d\zeta}\epsilon_{2n+1,1}(\mu,\alt,\zeta)\Bigl|_{\zeta=\zeta_0}\ll
\frac{\left|\tilde{I}_{i\tau}(\mu\sqrt{\zeta_0})\right|+\mu\left|\tilde{I}'_{i\tau}(\mu\sqrt{\zeta_0})\right|}{\mu^{2n+1}}.
\end{equation}
Therefore,
\begin{equation}
\mathcal{W}\left(\epsilon_{2n+1,1},W_{2n+1,2}\right)(\mu,\zeta_0)\ll
\frac{\left|\tilde{I}'_{i\tau}(\mu\sqrt{\zeta_0})K_{i\tau}(\mu\sqrt{\zeta_0})\right|+
\left|\tilde{I}_{i\tau}(\mu\sqrt{\zeta_0})K'_{i\tau}(\mu\sqrt{\zeta_0})\right|}{\mu^{2n}},
\end{equation}
and applying \cite[(4.25),(4.28)]{Dun90Bessel} one has
\begin{equation}\label{Wronskian e1W2}
\mathcal{W}\left(\epsilon_{2n+1,1},W_{2n+1,2}\right)(\mu,\zeta_0)=
O\left(\frac{e^{-\pi\tau}}{\mu^{2n+1}}\right).
\end{equation}
Substituting \eqref{CeCo-Ceco}, \eqref{Wronskian W1W2 asympt} to \eqref{e to W1W2} and using \eqref{Wronskian e1W2} we infer
\begin{multline}\label{e to W1W2 2}
e^{i\mu}_{-1/2+i\tau}(x)=\frac{1-\alt^2}{2\pi}e^{\pi\tau}\left(\frac{\zeta-\alt^2}{\zeta^2(1-\alt^2-\alt^2x^2)}\right)^{1/4}
\left(1+\sum_{s=1}^{n}\frac{D_s(\alt,\zeta_0)}{\mu^{2s}}+O\left(\frac{1}{\mu^{2n+1}}\right)\right)\times\\
\left(W_{2n+1,1}(\mu,\zeta)C^o_{2n+1,2}(\mu)-W_{2n+1,2}(\mu,\zeta)C^o_{2n+1,1}(\mu)\right).
\end{multline}
Applying \eqref{e to W1W2 2} and  \eqref{2F1 to e} we obtain
\begin{multline}\label{2f1 1/4 1/4 1/2 asympt1}
\HyG\left(1/4+ir(1-\alt),1/4-ir(1+\alt),1/2;y \right)=
\frac{(1-\alt^2)}{2\pi(1-y)^{-ir\alt}}e^{2\pi r\alt}
\left(\frac{\zeta-\alt^2}{\zeta^2(1-\alt^2-y)}\right)^{1/4}\times\\
\left(W_{2n+1,1}(2r,\zeta)C^o_{2n+1,2}(2r)-W_{2n+1,2}(2r,\zeta)C^o_{2n+1,1}(2r)\right)
\left(1+\sum_{s=1}^{n}\frac{D_s(\alt,\zeta_0)}{r^{2s}}+O\left(\frac{1}{r^{2n+1}}\right)\right).
\end{multline}

It is left to derive an asymptotic expression for $C^o_{2n+1,j}(2r)$.  To this end, we use \cite[(2.8)]{Dun90} and \cite[(4.28)]{Dun90Bessel}, proving that
\begin{equation}\label{I asympt}
\tilde{I}_{i\nu}(z)=\frac{\sqrt{2\pi}}{(z^2-\nu^2)^{1/4}} e^{-\pi\nu/2}e^{\nu\eta(z/\nu)}\left(1+\sum_{j=1}^{n-1}\frac{V_j(q)}{\nu^j}+O(\nu^{-n})\right),
\end{equation}
\begin{equation}\label{Ideriv asympt}
\tilde{I}'_{i\nu}(z)=\frac{(z^2-\nu^2)^{1/4}\sqrt{2\pi}}{z}e^{-\pi\nu/2}e^{\nu\eta(z/\nu)}\left(1+\sum_{j=1}^{n-1}\frac{\tilde{V}_j(z,q)}{\nu^j}+O(\nu^{-n})\right),
\end{equation}
where $$q=(z^2\nu^{-2}-1)^{-1/2}$$ and
\begin{equation}\label{eta def}
\eta(x)=\sqrt{x^2-1}-\arctan\sqrt{x^2-1}.
\end{equation}
Substituting \eqref{I asympt} and \eqref{Ideriv asympt} to \eqref{W1 def}, we infer
\begin{equation}\label{W1zeta0}
W_{2n+1,1}(\mu,\zeta_0)=
\frac{\sqrt{2\pi\zeta_0}}{(\zeta_0-\alt^2)^{1/4}\sqrt{\mu}} e^{-\pi\tau/2}e^{\tau\eta(\zeta_0^{1/2}/\alt)}\left(1+\sum_{j=1}^{n-1}\frac{a_j(\alt)}{\tau^j}+O(\tau^{-n})\right),
\end{equation}
\begin{equation}\label{W1derivzeta0}
W'_{2n+1,1}(\mu,\zeta_0)=
\frac{\sqrt{2\pi}(\zeta_0-\alt^2)^{1/4}\sqrt{\mu}}{2\sqrt{\zeta_0}}e^{-\pi\tau/2}e^{\tau\eta(\zeta_0^{1/2}/\alt)}\left(1+\sum_{j=1}^{n-1}\frac{b_j(\alt)}{\tau^j}+O(\tau^{-n})\right).
\end{equation}
Applying \eqref{zeta_0} one has
\begin{equation}\label{eta zeta0}
\eta\left(\frac{\zeta_0^{1/2}}{\alt}\right)=\frac{\pi(1-\alt)}{2\alt}.
\end{equation}
As a result, substituting \eqref{W1zeta0},  \eqref{W1derivzeta0} to \eqref{Co def}, we conclude that
\begin{equation}\label{Co1 asympt}
C^o_{2n+1,1}(2r)=
\frac{-\sqrt{4\pi r}}{(1-\alt^2)^{3/4}}e^{\pi r(1-2\alt)}\left(1+\sum_{j=1}^{n-1}\frac{c_j(\alt)}{(r\alt)^j}+O((r\alt)^{-n})\right).
\end{equation}
Using \cite[(4.25)]{Dun90Bessel} it is possible to obtain the asymptotic expansions (see also \cite{Balogh}):
\begin{equation}\label{K asympt}
K_{i\nu}(z)=\frac{\sqrt{\pi/2}}{(z^2-\nu^2)^{1/4}} e^{-\pi\nu/2}e^{-\nu\eta(z/\nu)}\left(1+\sum_{j=1}^{n-1}\frac{(-1)^jV_j(q)}{\nu^j}+O(\nu^{-n})\right),
\end{equation}
\begin{equation}\label{Kderiv asympt}
K'_{i\nu}(z)=\frac{-(z^2-\nu^2)^{1/4}\sqrt{\pi/2}}{z}e^{-\pi\nu/2}e^{-\nu\eta(z/\nu)}\left(1+\sum_{j=1}^{n-1}\frac{\breve{V}_j(z,q)}{\nu^j}+O(\nu^{-n})\right).
\end{equation}
Combining \eqref{K asympt}, \eqref{Kderiv asympt}, \eqref{W2 def}  and \eqref{Co def}, we show that
\begin{equation}\label{Co2 asympt}
C^o_{2n+1,2}(2r)=
\frac{\sqrt{\pi r}}{(1-\alt^2)^{3/4}}e^{-\pi r}\left(1+\sum_{j=1}^{n-1}\frac{d_j(\alt)}{(r\alt)^j}+O((r\alt)^{-n})\right).
\end{equation}
Finally, substituting \eqref{Co1 asympt}, \eqref{Co2 asympt} to \eqref{2f1 1/4 1/4 1/2 asympt1}, we prove \eqref{2F1 main asymptI}.

\section{Proof of Corollaries }
Using  Stirling's formula:
\begin{multline}\label{Stirling2}
\Gamma(\sigma+it)=\sqrt{2\pi}|t|^{\sigma-1/2}\exp(-\pi|t|/2)
\exp\left(i\left(t\log|t|-t+\frac{\pi t(\sigma-1/2)}{2|t|}\right)\right)\\\times
\left(1+\sum_{j=1}^{N-1}a_j/t^j+O(|t|^{-N})\right),
\end{multline}
which holds for $|t|\rightarrow\infty$ and a  fixed $\sigma$, one has
\begin{multline}\label{Gamma1/4Gamma1/4Stirling}
\frac{\Gamma(1/4+ir(1-\alt))\Gamma(1/4-ir(1+\alt))}{\Gamma(1/2)}=
\frac{2\sqrt{\pi}e^{-\pi r}}{(1-\alt^2)^{1/4}\sqrt{r}}\\\times
e^{ir\left(-2\alt\log r+(1-\alt)\log(1-\alt)-(1+\alt)\log(1+\alt)+2\alt\right)}
\left(1+\sum_{j=1}^{N-1}a_j/r^j+O(|r|^{-N})\right).
\end{multline}
Consider the case  $0<y<1-\alt^2-\frac{\alt^{4/3}}{r^{2/3-\delta}}$.
First, we would like to show that the contribution of $\tilde{I}_{2ir\alt}(2r\sqrt{\zeta})$ to \eqref{2F1 main asymptI} is negligible.
Since the function $\eta(x)$ (see \eqref{eta def}) increases and the identity  \eqref{eta zeta0} holds, we conclude that
\begin{equation}
\eta\left(\frac{\sqrt{\zeta}}{\alt}\right)\le\frac{\pi(1-\alt)}{2\alt}.
\end{equation}
Therefore, it follows from  \eqref{I asympt} and \eqref{K asympt}  that
\begin{equation}\label{I<K}
e^{-2\pi r+2\pi r\alt}\tilde{I}_{2ir\alt}(2r\sqrt{\zeta})\ll K_{2ir\alt}(2r\sqrt{\zeta}).
\end{equation}
Note that the left-hand side of \eqref{I<K} is approximately equal to the right-hand side only if $\zeta$ is close to $\zeta_0$, and in this case both parts
are of size $e^{-\pi r}$, and thus the part of \eqref{2F1 main asymptI} with $\tilde{I}_{2ir\alt}(2r\sqrt{\zeta})$  is either smaller than the part with
$K_{2ir\alt}(2r\sqrt{\zeta})$ or satisfies the estimate \eqref{2f1 1/4 1/4 1/2 asympt y<1-alt^2}. Therefore, it follows from \eqref{2F1 main asymptI} and \eqref{Gamma1/4Gamma1/4Stirling} that for $0<y<1-\alt^2-\frac{\alt^{4/3}}{r^{2/3-\delta}}$ the estimate holds:
\begin{equation}\label{2F1 est1 y<1-alt^2}
\Fs(r,\alt,y)\ll \left(\frac{\zeta-\alt^2}{1-\alt^2-y}\right)^{1/4}K_{2ir\alt}(2r\sqrt{\zeta}).
\end{equation}
In the case when $1-\alt^2-\frac{\alt^{4/3}}{r^{2/3-\delta}}<y<1$, the variable $\zeta$ is small, and thus $\tilde{I}_{2ir\alt}(2r\sqrt{\zeta})$ is absorbed by the error term $e^{\pi r}K_{2ir\alt}(2r\sqrt{\zeta})O(r^{-n})$ in  \eqref{2F1 main asymptI}.
As a result, the leading coefficient in the asymptotic expansion is given by
\begin{equation}\label{2f1 1/4 1/4 1/2 asympt2}
\Fs(r,\alt,y)\sim
e^{ir\Lt(\alt,y)}
\left(\frac{\zeta-\alt^2}{1-\alt^2-y}\right)^{1/4}
2K_{2ir\alt}(2r\sqrt{\zeta}),
\end{equation}
where $\Lt(\alt,y)$ is defined by \eqref{2F1 main asymptII Lt func}.
To deal with the $K$-Bessel function in  \eqref{2f1 1/4 1/4 1/2 asympt2} we apply its  well known (see \cite[(4.8)]{Dun90Bessel} or \cite{Balogh}) approximation in terms of Airy functions. The main term  is
\begin{equation}\label{K asympt1}
K_{i\nu}(\nu z)\sim\frac{\pi e^{-\pi\nu/2}}{\nu^{1/3}}
\left(\frac{4\hat{\zeta}}{1-z^2}\right)^{1/4}Ai(-\nu^{2/3}\hat{\zeta}),
\end{equation}
where for $0\le z\le1$
\begin{equation}\label{hatzeta z<1}
\frac{2}{3}\hat{\zeta}^{3/2}=\log\frac{1+\sqrt{1-z^2}}{z}-\sqrt{1-z^2},
\end{equation}
and for $z>1$
\begin{equation}\label{hatzeta z>1}
\frac{2}{3}(-\hat{\zeta})^{3/2}=\sqrt{z^2-1}-\arccos\left(\frac{1}{z}\right).
\end{equation}
Consequently, applying \eqref{K asympt1} with $z=\sqrt{\zeta}/\alt$ one has
\begin{equation}\label{2f1 1/4 1/4 1/2 asympt4}
\Fs(r,\alt,y)\sim
e^{ir\Lt(\alt,y)}
\frac{ 2^{3/2}\pi e^{-\pi r\alt}}{(2r\alt)^{1/3}}\left(\frac{\alt^2\hat{\zeta}(y)}{y-1+\alt^2}\right)^{1/4}Ai(-(2r\alt)^{2/3}\hat{\zeta}(y)).
\end{equation}
It follows from \eqref{zeta>alt^2} and \eqref{hatzeta z>1} that for $y<1-\alt^2$
\begin{equation}\label{hat zeta y<1-alt^2}
\frac{2\alt}{3}(-\hat{\zeta})^{3/2}(y)=
\arccos\frac{\sqrt{y}}{\sqrt{1-\alt^2}}-\frac{\alt}{2}\arccos\left(\frac{y(1+\alt^2)-(1-\alt^2)}{(1-\alt^2)(1-y)}\right)=\Aa_0(\alt,y),
\end{equation}
where $\Aa_0(\alt,y)$ is defined by \eqref{Aa0 def}.
Using \eqref{zeta<alt^2} and \eqref{hatzeta z<1} we show that for $1-\alt^2\le y<1$
\begin{equation}\label{hat zeta y>1-alt^2}
\frac{2\alt}{3}\hat{\zeta}(y)^{3/2}=
\log\frac{\sqrt{y}-\sqrt{y-1+\alt^2}}{\sqrt{1-\alt^2}}-
\frac{\alt}{2}\log\left(\frac{2\alt\left(\alt-\sqrt{y(y-1+\alt^2)}\right)}{(1-\alt^2)(1-y)}-\frac{1+\alt^2}{1-\alt^2}\right)=\Aa_1(\alt,y),
\end{equation}
where $\Aa_1(\alt,y)$ is defined by \eqref{Aa1 def}. Thus the Corollary \ref{cor: y=1-alt^2} is proved.

For the Airy function one can use asymptotic formulas \cite[(9.7.5),(9.7.9)]{HMF} provided that its argument is bigger than $r^{\delta}$. For $y<1-\alt^2$, $r\Aa_0(y)\gg r^{\delta}$ one has
\begin{equation}\label{Ai y<1-alt^2}
Ai(-(2r\alt)^{2/3}\hat{\zeta}(y))=Ai((3r\Aa_0(\alt,y))^{2/3})\sim \frac{e^{-2r\Aa_0(\alt,y)}}{2\sqrt{\pi}(3r\Aa_0(\alt,y))^{1/6}},
\end{equation}
and for $y>1-\alt^2$, $r\Aa_1(y)\gg r^{\delta}$  one has
\begin{equation}\label{Ai y>1-alt^2}
Ai(-(2r\alt)^{2/3}\hat{\zeta}(y))=Ai(-(3r\Aa_1(\alt,y))^{2/3})\sim \frac{\cos\left(2r\Aa_1(\alt,y)-\pi/4\right)}{\sqrt{\pi}(3r\Aa_1(\alt,y))^{1/6}}.
\end{equation}
Since $\Aa_0(\alt,1-\alt^2)=\Aa_1(\alt,1-\alt^2)=0$, we express these functions in the Taylor series at the point $y=1-\alt^2$. Accordingly,
\begin{multline}\label{Aa0 Taylor}
\Aa_0(\alt,y)=
\frac{\pi(1-\alt)}{2}-\arctan\frac{\sqrt{y}}{\sqrt{1-\alt^2-y}}+\alt\arctan\frac{\alt\sqrt{y}}{\sqrt{1-\alt^2-y}}=\\=
\frac{(1-\alt^2)(1-\alt^2-y)^{3/2}}{3\alt^2y^{3/2}}\left(1+O\left(\frac{1-\alt^2-y}{\alt^2}\right)\right),
\end{multline}
\begin{multline}\label{Aa1 Taylor}
\Aa_1(\alt,y)=
\alt\log\left(\alt\sqrt{y}+\sqrt{y-1+\alt^2}\right)-\log\left(\sqrt{y}+\sqrt{y-1+\alt^2}\right)-\frac{\alt}{2}\log(1-y)+
\frac{1-\alt}{2}\log(1-\alt^2)=\\=
\frac{(y-1+\alt^2)^{3/2}}{3\alt^2\sqrt{1-\alt^2}}\left(1+O\left(\frac{y-1+\alt^2}{\alt^2}\right)\right).
\end{multline}
Therefore, for
\begin{equation}\label{y-1+alt^2  condition}
|y-1+\alt^2|\gg\frac{\alt^{4/3}}{r^{2/3-\delta}}
\end{equation}
one has $r\Aa_j(y)\gg r^{3\delta/2}$, and consequently, it is possible to substitute \eqref{Ai y<1-alt^2}, \eqref{Ai y>1-alt^2} to \eqref{2f1 1/4 1/4 1/2 asympt4}, obtaining
\begin{equation}\label{2f1 1/4 1/4 1/2 asympt5 y<1-alt^2}
\Fs(r,\alt,y)\sim
(2/3)^{1/6}\sqrt{\pi}
e^{ir\Lt(\alt,y)}
\frac{ e^{-\pi r\alt}}{(r\alt)^{1/3}}\left(\frac{\alt^2\hat{\zeta}(y)}{y-1+\alt^2}\right)^{1/4}
\frac{e^{-2r\Aa_0(\alt,y)}}{(r\Aa_0(\alt,y))^{1/6}}
\end{equation}
for $0<y<1-\alt^2-\frac{\alt^{4/3}}{r^{2/3-\delta}}$ and
\begin{equation}\label{2f1 1/4 1/4 1/2 asympt5 y>1-alt^2}
\Fs(r,\alt,y)\sim
2(2/3)^{1/6}\sqrt{\pi}
e^{ir\Lt(\alt,y)}
\frac{ e^{-\pi r\alt}}{(r\alt)^{1/3}}\left(\frac{\alt^2\hat{\zeta}(y)}{y-1+\alt^2}\right)^{1/4}
\frac{\cos\left(2r\Aa_1(\alt,y)-\pi/4\right)}{(r\Aa_1(\alt,y))^{1/6}}
\end{equation}
$1-\alt^2+\frac{\alt^{4/3}}{r^{2/3-\delta}}<y<1.$
Finally, substituting  \eqref{hat zeta y<1-alt^2} and \eqref{hat zeta y>1-alt^2} to \eqref{2f1 1/4 1/4 1/2 asympt5 y<1-alt^2} and \eqref{2f1 1/4 1/4 1/2 asympt5 y>1-alt^2}, we prove \eqref{2f1 1/4 1/4 1/2 asympt y<1-alt^2} and \eqref{2f1 1/4 1/4 1/2 asympt y>1-alt^2}.
\section{Numerical examples}
In this section, we provide some numerical examples related to Corollaries \ref{cor: y>1-alt^2} and \ref{cor: y=1-alt^2}.
The results are given in Tables \ref{case1}-\ref{case6}.

\begin{center}
 \begin{tabular}{ |c|c|c| }
\hline
Function & Approximation & $|$ Relative error $|$ \\ \hline
 &  &  \\
$\Fs(r,\alt,y)$  & $-6.008705138*10^{-16}-1.461019048*10^{-15}i$ &  \\
&  &  \\
rhs \eqref{2f1 1/4 1/4 1/2 asympt y>1-alt^2} & $-5.595045762*10^{-16}-1.360356225*10^{-15}i $ & $0.068890994$  \\
&  &  \\
rhs \eqref{2f1 1/4 1/4 1/2 asympt y=1-alt^2} & $-5.763136445*10^{-16}-1.401225097*10^{-15}i $ & $0.0409179003$  \\
&  &  \\ \hline
\end{tabular}
\captionof{table}{The case $r=100$, $\alpha=0.1$, $y=0.9999$.}\label{case1}
\end{center}


\begin{center}
 \begin{tabular}{ |c|c|c| }
\hline
Function & Approximation & $|$ Relative error $|$ \\ \hline
 &  &  \\
$\Fs(r,\alt,y)$  & $1.854052580*10^{-14}-2.566435037*10^{-14}i$ &  \\
&  &  \\
rhs \eqref{2f1 1/4 1/4 1/2 asympt y>1-alt^2} & $2.246324758*10^{-14}-3.109567770*10^{-14}i $ & $0.211610836$  \\
&  &  \\
rhs \eqref{2f1 1/4 1/4 1/2 asympt y=1-alt^2} & $1.854248101*10^{-14}-2.566819473*10^{-14}i $ & $0.000136225$  \\
&  &  \\ \hline
\end{tabular}
\captionof{table}{The case $r=100$, $\alpha=0.1$, $y=0.991$.}\label{case2}
\end{center}


\begin{center}
 \begin{tabular}{ |c|c|c| }
\hline
Function & Approximation & $|$ Relative error $|$ \\ \hline
 &  &  \\
$\Fs(r,\alt,y)$  & $-2.795222815*10^{-14}+9.161763311*10^{-15}i$ &  \\
&  &  \\
rhs \eqref{2f1 1/4 1/4 1/2 asympt y>1-alt^2} & $-2.755422703*10^{-14}+9.031954601*10^{-15}i $ & $0.014231835$  \\
&  &  \\
rhs \eqref{2f1 1/4 1/4 1/2 asympt y=1-alt^2} & $-2.794277351*10^{-14}+9.159315612*10^{-15}i $ & $0.000332015$  \\
&  &  \\ \hline
\end{tabular}
\captionof{table}{The case $r=100$, $\alpha=0.1$, $y=0.993$.}\label{case3}
\end{center}



\begin{center}
 \begin{tabular}{ |c|c|c| }
\hline
Function & Approximation & $|$ Relative error $|$ \\ \hline
 &  &  \\
$\Fs(r,\alt,y)$  & $7.792063111*10^{-15}+1.395437691*10^{-14}i$ &  \\
&  &  \\
rhs \eqref{2f1 1/4 1/4 1/2 asympt y=1-alt^2} & $7.794616502*10^{-15}+1.395818219*10^{-14}i $ & $0.000285226$  \\
&  &  \\ \hline
\end{tabular}
\captionof{table}{The case $r=100$, $\alpha=0.1$, $y=0.989$.}\label{case4}
\end{center}


\begin{center}
 \begin{tabular}{ |c|c|c| }
\hline
Function & Approximation & $|$ Relative error $|$ \\ \hline
 &  &  \\
$\Fs(r,\alt,y)$  & $0.001824209+0.002110461i$ &  \\
&  &  \\
rhs \eqref{2f1 1/4 1/4 1/2 asympt y>1-alt^2} & $0.001725718+0.001996498i $ & $0.053995710$  \\
&  &  \\
rhs \eqref{2f1 1/4 1/4 1/2 asympt y=1-alt^2} & $0.001807240+0.002090810i $ & $0.009307311$  \\
&  &  \\ \hline
\end{tabular}
\captionof{table}{The case $r=100$, $\alpha=0.02$, $y=0.9999$.}\label{case5}
\end{center}



\begin{center}
 \begin{tabular}{ |c|c|c| }
\hline
Function & Approximation & $|$ Relative error $|$ \\ \hline
 &  &  \\
$\Fs(r,\alt,y)$  & $-0.004422131+0.000382789i$ &  \\
&  &  \\
rhs \eqref{2f1 1/4 1/4 1/2 asympt y>1-alt^2} & $-0.005487073+0.000474996i $ & $0.240821678$  \\
&  &  \\
rhs \eqref{2f1 1/4 1/4 1/2 asympt y=1-alt^2} & $-0.004427232+0.000383249i $ & $0.001154018$  \\
&  &  \\ \hline
\end{tabular}
\captionof{table}{The case $r=100$, $\alpha=0.02$, $y=0.9997$.}\label{case6}
\end{center}


\end{document}